\documentclass[11pt]{amsart}
\usepackage{amsmath}
\usepackage{amssymb}
\usepackage{amscd}
\usepackage{color}

\def\NZQ{\mathbb}               
\def\NN{{\NZQ N}}

\def\ZZ{{\NZQ Z}}
\def\RR{{\NZQ R}}

\newtheorem{Theorem}{Theorem}[section]

\newtheorem{Proposition}[Theorem]{Proposition}
\newtheorem{Remark}[Theorem]{Remark}

\let\epsilon\varepsilon
\let\phi=\varphi
\let\kappa=\varkappa

\begin{document}

\title{A counterexample to strong local monomialization in a tower of two independent defect Artin-Schreier extensions}
\author{Steven Dale Cutkosky}
\thanks{Partially supported by NSF grant DMS-2054394}

\address{Steven Dale Cutkosky, Department of Mathematics,
University of Missouri, Columbia, MO 65211, USA}
\email{cutkoskys@missouri.edu}


\begin{abstract} We give an example of an extension of two dimensional regular local rings in a tower of two independent defect Artin-Schreier extensions for which strong local monomialization does not hold. 
\end{abstract}

\maketitle 

\section{Introduction}

In characteristic zero, there is a very nice local form for morphisms, called local monomialization. This result is a little stronger than what comes immediately from the assumption that toroidalization is possible.  If $R\rightarrow S$ is an extension of local rings such that the maximal ideal of $S$ contracts to the maximal ideal of $R$ then we say that $S$ dominates $R$. If $S$ is dominated by the valuation ring $\mathcal O_{\omega}$ of a valuation $\omega$ we say that $\omega$ dominates $S$.

\begin{Theorem}\label{locmon}(local monomialization)(\cite{C}, \cite{C5}) Suppose that $k$ is a field of characteristic zero and $R\rightarrow S$ is an extension of regular local rings such that $R$ and $S$ are essentially of finite type over $k$ and $\omega$ is a valuation of the quotient field of $S$ which dominates $S$ and $S$ dominates $R$. Then there is a 
commutative diagram
$$
\begin{array}{ccc}
R_1&\rightarrow & S_1\\
\uparrow&&\uparrow\\
R&\rightarrow &S
\end{array}
$$
such that  $\omega$ dominates $S_1$, $S_1$ dominates $R_1$ and the vertical arrows are products of monoidal transforms; that is, these arrows are factored by the local rings of blowups of prime ideals  whose quotients are regular local rings. In particular, $R_1$ and $S_1$ are regular local rings. Further, $R_1\rightarrow S_1$ has a locally monomial form; that is, there exist regular parameters $u_1,\ldots,u_m$ in $R_1$ and $x_1,\ldots,x_n$ in $S_1$, an $m\times n$ matrix $A=(a_{ij})$ with integral coefficients such that 
$\mbox{rank}(A)=m$ and units $\delta_i\in S_1$ such that
$$
u_i=\delta_i\prod_{j=1}^nx_j^{a_{ij}}
$$
for $1\le i\le m$.
\end{Theorem}
The difficulty in the proof is to obtain the condition that $\mbox{rank}(A)=m$. To do this, it is necessary to blow up above both $R$ and $S$.

In the case when the extension of quotient fields $K\rightarrow L$ of the extension $R\rightarrow S$ is a finite extension and $k$ has characteristic zero, it is possible to find a local monomialization such that the structure of the matrix of coefficients recovers classical invariants of the extension of valuations in $K\rightarrow L$, and this form holds stably along suitable sequences of birational morphisms which generate the respective valuation rings. This form is called strong local uniformization. It is established for rank 1 valuations in \cite{C} and for general valuations in \cite{CP}. The case which has the simplest form and will be of interest to us in this paper is when the valuation has rational rank 1. In this case, if $R_1\rightarrow S_1$ is a strong local monomialization, then there exist regular parameters $u_1,\ldots, u_m$ in $R_1$ and $v_1,\ldots,v_m$ in $S_1$ , a positive integer $a$ and a unit $\delta\in S_1$ such that
\begin{equation}\label{N20}
u_1=\delta v_1^a, u_2=v_2,\ldots,u_m=v_m.
\end{equation}

The stable forms of mappings in positive characteristic and dimension $\ge 2$ are much more complicated. For instance, local monomialization does not always hold. An example is given in \cite{C2} where $R\rightarrow S$ are local rings of points on nonsingular algebraic surfaces over an algebraically closed field $k$ of positive characteristic $p$ and $k(X)\rightarrow k(Y)$ is finite and separable.

The obstruction to local monomialization is the defect.   The defect $\delta(\omega/\nu)$, which is a power of the residue characteristic $p$ of $\mathcal O_{\omega}$,  is defined and its basic properties developed in \cite[Chapter VI, Section 11]{ZS2}, \cite{K0}, \cite[Section 7.1]{CP}. The defect is discussed in Subsection \ref{Subsecnot}. We have the following theorem, showing that the defect is the only obstruction to strong local monomialization for maps of surfaces.

\begin{Theorem}\label{nondef}(\cite[Theorem 7.35]{CP}) Suppose that $K\rightarrow L$ is a finite, separable extension of algebraic function fields over an algebraically closed field $k$ of characteristic $p>0$, $R\rightarrow S$ is an extension of local domains such that $R$ and $S$ are essentially of finite type over $k$ and the quotient fields of $R$ and $S$ are $K$ and $L$ respectively such that $S$ dominates $R$. Suppose that $\omega$ is valuation of $L$ which dominates $S$. Let $\nu$ be the restriction of $\omega$ to $K$. Suppose that the extension is defectless ($\delta(\omega/\nu)=1$). Then the conclusions of Theorem \ref{locmon} hold. In particular, $R\rightarrow S$ has a local monomialization (and a strong local monomialization) along $\omega$. 
\end{Theorem}

Suppose that 
$K\rightarrow L$ is a Galois extension of fields of characteristic $p>0$  and $\omega$ is a valuation of $L$, $\nu$ is the restriction of $\omega$ to $K$. Then there is a classical tower of fields (\cite[page 171]{End})
$$
K\rightarrow K^s\rightarrow K^i\rightarrow K^v\rightarrow L.
$$
where $K^s$ is the splitting field, $K^i$ is the inertia field, $K^v$ is the ramification field and the extension $K\rightarrow K^v$ has no defect. Thus the essential difficulty comes from the extension from  $K^v$ to $L$ which could have defect. The extension $K^v\rightarrow L$ is a tower of Artin-Schreier extensions, so the Artin-Schreier extension is of fundamental importance in this theory.  

Kuhlmann has extensively studied defect in Artin-Schreier extensions in \cite{Ku}. He separated these extensions into  dependent and independent defect Artin-Schreier extensions. This definition is reproduced in Subsection \ref{Galois}.
Kuhlmann  also defined an invariant called the distance to distinguish the natures of Artin-Schreier extensions. This definition is  given in Subsections \ref{SecDist} and \ref{Galois}.

We now specialize to the case of  a finite separable  extension $K \rightarrow L$ of two   dimensional algebraic function fields  over an algebraically closed field $k$ of characteristic $p>0$, and suppose that $\omega$ is a valuation of $L$ which is trivial on $k$ and $\nu$ is the restriction of $\omega$ to $K$.  If $L/K$ has defect then $\omega$ must have  rational rank 1 and be nondiscrete. We will assume that $\omega$ has rational rank 1 and is nondiscrete for the remainder of the introduction.

 With these restrictions, the distance $\delta$ of an Artin-Schreier extension 
is $\le$ $0^{-}$ when the extension has  defect. If it is a defect extension with $\delta=0^{-}$  then it is  an independent defect extension. If it is a defect extension and the distance is less than $0^{-}$ then the extension is a dependent defect extension.

A quadratic transform along a valuation is the center of the valuation at the blow up of a maximal ideal of a regular local ring. There is the sequence of quadratic transforms along $\nu$ and $\omega$
\begin{equation}\label{In2}
R\rightarrow R_1\rightarrow R_2\rightarrow \cdots\mbox{ and }S\rightarrow S_1\rightarrow S_2\rightarrow\cdots.
\end{equation}
We have that $\cup_{i=1}^{\infty}R_i=\mathcal O_{\nu}$, the valuation ring of $\nu$, and $\cup_{i=1}^{\infty}S_i=\mathcal O_{\omega}$, the valuation ring of $\omega$.
 These sequences can be factored by standard quadratic transform sequences (defined in Section \ref{Sec2Ext}).  It is shown in \cite{CP} that given positive integers $r_0$ and $s_0$, there exists $r\ge r_0$ and $s\ge s_0$ such that $R_r\rightarrow S_s$ has the following  form:
\begin{equation}\label{In1}
u=\delta x^a, v=x^b(y^d\gamma+x\Omega)
\end{equation}
where $u,v$ are regular parameters in $R_r$, $x,y$ are regular parameters in $S_s$, $\gamma$ and $\tau$ are units in $S_s$,
$\Omega\in S_s$, $a$ and $d$ are positive integers and $b$ is a non negative integer. If we choose $r_0$  sufficiently large, then
we have  that the complexity $ad$ of the extension $R_r\rightarrow S_s$  is a constant which depends on the extension of valuations, which we call the stable complexity of (\ref{In2}). When $R_r\rightarrow S_s$ has this stable complexity, we  call the forms (\ref{In1}) stable forms.

  The strongly monomial form is the case when $b=0$ and $d=1$; that is, after making a change of variables in $y$,
  $$
  u=\delta x^a, v=y.
  $$

As we observed earlier (Theorem \ref{nondef}) if the extension $K\rightarrow L$ has no defect, then the stable form is the strongly monomial form. If there is defect, then it is possible for the $a$ and $d$ in stable forms along a valuation to vary wildly, even though their product $ad$ is fixed by the extension, as shown in \cite[Theorem 5.4]{C3}.

An example is constructed in \cite{CP}, showing failure of strong local monomialization.    It is a tower of two defect Artin-Schreier extensions, each  of the type of \cite[Theorem 5.4]{C3} referred to above. The first extension is of type 1 for even integers and of type 2 for odd integers. The second extension is of type 2 for even integers and of type 1 for odd integers. The composite gives a sequence of extensions of regular local  rings $R_i\rightarrow S_i$, where $R_i$ has regular parameters $u_i,v_i$ and $S_i$ has regular parameters $x_i,y_i$ such that the stable form is 
\begin{equation}\label{eqN25}
u_i=\gamma x_i^{p}, v_i=y_i^p\tau+x_i\Omega
\end{equation}
for all $i$.
Both of these  Artin-Schreier extensions are  dependent. This is calculated in \cite{EG} and in \cite[Section 6]{C3}.
In keeping with the philosophy that independent Artin-Schreier extensions are better behaved than dependent ones, this leads to the question of if  strong monomialization holds in towers of independent Artin-Schreier extensions.  However, this is not true as is shown in Theorem \ref{Example3} of this paper. In this theorem,  we  construct an example in   a tower of two independent defect extensions such that  strong local monomialization does not hold. 

Suppose that  $K\rightarrow L$ is a finite extension of  fields of positive characteristic and $\omega$ is a valuation of $L$  with restriction $\nu$ to $K$. It is known that there is no defect in the extension if and only if there is a finite generating sequence in $L$ for the valuation $\omega$ over $K$ (\cite{Va}, \cite{NS}). The calculation of generating sequences for extensions of Noetherian local rings which are dominated by a valuation is extremely difficult. This has been accomplished for two  dimensional  regular local rings  in \cite{Sp} and \cite{CV} and for many hypersurface singularities above a regular local ring of arbitrary dimension in \cite{CMT}.

The nature of a generating sequence in an extension of  $S$ over $R$ determines the nature of the mappings in the stable forms. It is shown in \cite[Theorem 1]{C8} that if $R\rightarrow S$ is an extension of two dimensional excellent regular local rings whose quotient fields give a finite extension $K\rightarrow L$ and $\omega$ is a valuation of $L$ which dominates $S$ then the extension is without defect if and only if  there exist sequences of quadratic transform $R\rightarrow R_1$ and $S\rightarrow S_1$ along $\nu$ 
such that $\omega$ has  a finite generating sequence in $S_1$ over $R_1$.
This shows us that we can expect good stable forms (as do hold by Theorem \ref{nondef}) if there is no defect, but not otherwise.

\section{Preliminaries}\label{SecPre}
\subsection{Some notation}\label{Subsecnot}  Let $K$ be a field with a valuation $\nu$.  The valuation ring of $\nu$ will be donoted by $\mathcal O_{\nu}$, $\nu K$ will denote the value group of $\nu$ and $K\nu$ will denote the residue field of $\mathcal O_{\nu}$. 

The maximal ideal of a local ring $A$ will be denoted by $m_A$. If $A\rightarrow B$ is an extension (inclusion) of local rings such that $m_B\cap A=m_A$ we will say that $B$ dominates $A$. If a valuation ring $\mathcal O_{\nu}$ dominates $A$ we will say that the valuation $\nu$ dominates $A$.

Suppose that $K$ is an algebraic function field over a field $k$. An algebraic local ring $A$ of $K$ is a local domain which is a localization of a finite type $k$-algebra whose quotient field is $K$. A $k$-valuation of $K$ is a valuation of $K$ which is trivial on $k$.

Suppose that $K\rightarrow L$ is a finite algebraic extension of fields, $\nu$ is a valuation of $K$  and $\omega$ is an extension of $\nu$ to $L$. Then the reduced ramification index of the extension is $e(\omega/\nu)=[\omega L:\nu K]$ and the residue degree of the extension is $f(\omega/\nu)=[L\omega:K\nu]$.

 The defect $\delta(\omega/\nu)$, which is a power of the residue characteristic $p$ of $\mathcal O_{\omega}$,  is defined and its basic properties developed in \cite[Chapter VI, Section 11]{ZS2}, \cite{K0} and \cite[Section 7.1]{CP}.
 In the case that $L$ is Galois over $K$, we have the formula 
 \begin{equation}\label{found8}
 [L:K]=e(\omega/\nu)f(\omega/\nu)\delta(\omega/\nu)g
 \end{equation}
 where $g$ is the number of  extensions of $\nu$ to $L$. In fact, 
 we have  the equation (c.f. \cite{Ku} or Section 7.1 \cite{CP})
$$
|G^s(\omega/\nu)|=e(\omega/\nu)f(\omega/\nu)\delta(\omega/\nu),
$$
where $G^s(\omega/\nu)$ is the decomposition group of $L/K$. 

If $K\rightarrow L$ is a finite Galois extension, then we will denote the Galois group of $L/K$ by $\mbox{Gal}(L/K)$.

\subsection{Initial and final segments and cuts}
We review some basic material about cuts in totally ordered sets from \cite{Ku}.
Let $(S,<)$ be a totally ordered set. An initial segment of $S$ is a subset $\Lambda$ of $S$ such that if $\alpha\in \Lambda$ and $\beta<\alpha$ then $\beta\in \Lambda$. A final segment of $S$ is a subset $\Lambda$ of $S$ such that if $\alpha\in \Lambda$ and $\beta>\alpha$ then $\beta\in \Lambda$. A cut in $S$ is a pair of sets $(\Lambda^L,\Lambda^R)$ such that $\Lambda^L$ is an initial segment of $S$ and $\Lambda^R$ is a final segment of $S$ satisfying $\Lambda^L\cup \Lambda^R=S$ and $\Lambda^L\cap \Lambda^R=\emptyset$.
If $\Lambda_1$ and $\Lambda_2$ are two cuts in $S$, write $\Lambda_1<\Lambda_2$ if $\Lambda_1^L\subsetneq \Lambda_2^{L}$.
Suppose that $S\subset T$ is an order preserving inclusion of ordered sets and $\Lambda=(\Lambda^L,\Lambda^R)$ is a cut in $S$. Then define the cut induced by $\Lambda=(\Lambda^L,\Lambda^R)$ in $T$ to be the cut 
$\Lambda\uparrow T=(\Lambda^L\uparrow T,T\setminus (\Lambda^L\uparrow T))$ where $\Lambda^L\uparrow T$ is the least initial segment of $T$ in which $\Lambda^L$ forms a cofinal subset. 

We embed $S$ in the set of all cuts of $S$ by sending $s\in S$ to 
$$
s^+=(\{t\in S\mid t\le s\},\{t\in S\mid t>s\}).
$$
we may identify $s$ with the cut $s^+$. Define
$$
s^{-}=(\{t\in S\mid t< s\},\{t\in S\mid t\ge s\}).
$$
Given a cut $\Lambda=(\Lambda^L,\Lambda^R)$, we define 
$-\Lambda=(-\Lambda^R,-\Lambda^L)$ where
$-\Lambda^L=\{-s\mid s\in \Lambda^L\}$ and $-\Lambda^R=\{-s\mid s\in \Lambda^R\}$.
We have that if $\Lambda_1$ and $\Lambda_2$ are cuts, then $\Lambda_1<\Lambda_2$ if and only if $-\Lambda_2<-\Lambda_1$.

Observe that for $s\in S$, $-s=-(s^+)=(-s)^-$ and $-(s^-)=(-s)^+=-s$.

\subsection{Distances}\label{SecDist}
Let $K\rightarrow L$ be an   extension of  fields and  $\omega$ be a valuation of $L$ with restriction $\nu$ to $K$. 
Let $\widetilde{\nu K}$ be the divisible hull of  $\nu K$. Suppose that $z\in L$. Then the distance of $z$ from $K$ is defined in \cite[Section 2.3]{Ku}  to be the cut $\mbox{dist}(z,K)$ of $\widetilde{\nu K}$ in which  the initial segment of $\mbox{dist}(z,K)$ is the least initial segment of
$\widetilde{\nu K}$ in which $\omega(z-K)$ is cofinal. That is, 
$$
\mbox{dist}(z,K)=(\Lambda^L(z,K),\Lambda^R(z,K))\uparrow\widetilde{\nu K}
$$
where 
$$
\Lambda^L(z,K)=\{\omega(z-c)\mid c\in K\mbox{ and }\omega(z-c)\in \nu K\}.
$$
The following notion of equivalence is defined in \cite[Section 2.3]{Ku}. If $y,z\in L$, then $z\sim_Ky$ if  
$\omega(z-y) >\mbox{dist}(z,K)$.

\subsection{Artin-Schreier extensions}\label{Galois}
Let $K\rightarrow L$ be an  Artin-Schreier extension of  fields of characteristic $p>0$ and  $\omega$ be a valuation of $L$ with restriction $\nu$ to $K$. The field $L$ is Galois over $K$ with Galois group $G\cong\ZZ_p$, where $p$ is the characteristic of $K$.

Let $\Theta\in L$ be an Artin-Schreier generator of $K$; that is, there is an expression
$$
\Theta^p-\Theta = a
$$
for some $a\in K$. We have that 
$$
\mbox{Gal}(L/K)\cong \ZZ_p=\{{\rm id}, \sigma_1,\ldots,\sigma_{p-1}\},
$$
where $\sigma_i(\Theta)=\Theta+i$.

Since $L$ is Galois over $K$, we have that $ge(\omega/\nu)f(\omega/\nu)\delta(\omega/\nu)=p$ where $g$ is the number of extensions of $\nu$ to $L$. So  we either have that $g=1$ or $g=p$.
If $g=1$, then $\omega$ is the unique extension of $\nu$ to $L$ and either $e(\omega/\nu)=p$ and $\delta(\omega/\nu)=1$ or $e(\omega/\nu)=1$ and $\delta(\omega/\nu)=p$. In particular, the extension is defect if and only if is an immediate extension ($e=f=1$) and $\omega$ is the unique extension of $\nu$ to $L$.


From now on in this subsection, suppose that $L$ is a defect extension of $K$. By \cite[Lemma 4.1]{Ku}, the distance 
$\delta=\mbox{dist}(\Theta,K)$ does not depend on the choice of Artin-Schreier generator $\Theta$, so $\delta$ can be called the distance of the Artin-Schreier extension. Since $L/K$ is an immediate extension, the set $\omega(\Theta - K)$ is an initial segment in $\nu K$ which  has no maximal element by \cite[Theorem 2.19]{Ku}.

 We have, since the extension is defect, that
\begin{equation}\label{N22}
\delta={\rm dist}(\Theta,K)\le 0^-
\end{equation}
by \cite[Corollary 2.30]{Ku}.

A defect Artin-Schreier extension $L$ is defined in \cite[Section 4]{Ku} to be a dependent defect Artin-Schreier extension if there exists an immediate purely inseparable extension $K(\eta)$ of $K$ of degree $p$ such that $\eta\sim_K\Theta$. Otherwise, $L/K$ is defined to be an independent defect Artin-Schreier defect extension. We have by \cite[Proposition 4.2]{Ku} that for a defect Artin- Schreier extension,
\begin{equation}\label{eqN23}
\mbox{$L/K$ is independent if and only if the distance $\delta={\rm dist}(\Theta,K)$  satisfies $\delta=p\delta$.}
\end{equation}



\subsection{Extensions of rank 1 valuations in an Artin-Schreier extension}\label{Rank1AS}
In this subsection, we suppose that $L$ is an Artin-Schreier extension of a field $K$ of characteristic $p$, $\omega$ is a  rank 1 valuation of $L$ and $\nu$ is the restriction of $\omega$ to $K$. We suppose that $L$ is a defect extension of $K$. To simplify notation, we suppose that we have an embedding of $\nu L$ in $\RR$. 
Since $L$ has defect over $K$ and $L$ is separable over $K$, $\nu L$ is nondiscrete by the corollary on page 287 of \cite{ZS1}, so that $\nu L$ is dense in $\RR$.



We define a cut 
in $\RR$ by extending the cut $\mbox{dist}(\Theta,K)$ in $\widetilde{\nu K}$ to a cut of $\RR$ by taking the initial segment of 
 the extended cut to be the least initial segment of $\RR$ in which  the cut $\mbox{dist}(\Theta,K)$ is confinal.
 This cut is then $\mbox{dist}(\Theta,K)\uparrow \RR$.  This cut is either $s$ or $s^{-}$ for some $s\in \RR$. If $L$ is a defect extension of $K$ then $\mbox{dist}(\Theta,K)\uparrow \RR=s^{-}$ where
$s$ is a non positive real number by \cite[Theorem 2.19]{Ku} and \cite[Corollary 2.30]{Ku}. 
We will set 
$\mbox{dist}(\omega/\nu)$ to be this real number  $s$, so that 
$$
\mbox{dist}(\Theta,K)\uparrow \RR=s^{-}=(\mbox{dist}(\omega/\nu))^{-}.
$$
The real number $\mbox{dist}(\omega/\nu)$  is well defined since it is independent of choice of Artin-Schreier generator of $L/K$ by Lemma 4.1 \cite{Ku}.

With the assumptions of this subsection, by (\ref{N22}) and (\ref{eqN23}), the distance $\delta=\mbox{dist}(\Theta,K)$ of an Artin-Schreier extension 
is $\le$ $0^{-}$ when the extension has defect. If it is a defect extension with distance equal to $0^{-}$ then it is an  independent defect extension. If it is a defect extensions and the distance is less than $0^{-}$ then the extension is a dependent defect extension. Thus if $L/K$ is a defect extension, we have that $\mbox{dist}(\omega/\nu)\le 0$ and the defect extension $L/K$ is independent if and only if $\mbox{dist}(\omega/\nu)= 0$.

\section{Calculations in  two dimensional Artin-Schreier Extensions}\label{Sec2Ext}\label{SecCalc}

Suppose that $M$ is a two dimensional algebraic function field over an algebraically closed field $k$ of characteristic $p>0$ and $\mu$ is a nondiscrete rational rank 1 valuation of $M$.
Suppose that $A$ is an algebraic regular local ring of $M$ such that $\mu$ dominates $A$. A quadratic transform of $A$ is an extension $A\rightarrow A_1$ where $A_1$ is a local ring of the blowup of the maximal ideal of $A$ such that $A_1$ dominates $A$ and $A_1$ has dimension two. A quadratic transform $A\rightarrow A_1$ is said to be along the valuation $\mu$ if $\mu$ dominates $A_1$.

Let
$$
A=A_0\rightarrow A_1\rightarrow A_2\rightarrow \cdots
$$
be the sequence of quadratic transforms along $\mu$. Then  the valuation ring $\mathcal O_{\mu}=\cup A_i$ (by \cite[Lemma 12]{Ab1}).  

Suppose that $K\rightarrow L$ is a finite extension of two dimensional algebraic function fields,
 $R$ is an algebraic  regular local ring of $K$ which is dominated by a regular algebraic local ring $S$ of $L$ such that $\dim R=\dim S=2$. Let $x,y$ be regular parameters in $S$ and $u,v$ be regular parameters in $R$. Then we can form the Jacobian ideal 
$$
J(S/R)=
(\frac{\partial u}{\partial x}\frac{\partial v}{\partial y}-\frac{\partial u}{\partial y}\frac{\partial v}{\partial x}).
$$
This ideal is independent of choice of regular parameters.

The following proposition is proven in  \cite{Pi}.

\begin{Proposition}\label{Prop1*} 
Suppose that $K\rightarrow L$ is an Artin-Schreier extension of two dimensional algebraic function fields over an algebraically closed field $k$ of characteristic $p>0$, $\omega$ is a rational rank 1 nondiscrete valuation of $L$ with restriction $\nu=\omega|K$. Further suppose that $A$ is an algebraic  local ring of $K$ and $B$ is an algebraic local ring of $L$ which is dominated by $\omega$ such that $B$ dominates $A$. Then there exists a commutative diagram of homomorphisms
$$
\begin{array}{ccc}
R&\rightarrow &S\\
\uparrow&&\uparrow \\
A&\rightarrow &B
\end{array}
$$
such that $R$ is a regular algebraic local ring of $K$ with regular parameters $u,v$, $S$ is a regular algebraic local ring of $L$ with regular parameters $x,y$ such that $S$ is dominated by $\omega$, $S$ dominates $R$,  $R\rightarrow S$ is quasi finite,
$J(S/R)= (x^{\overline c})$ for some non negative integer $\overline c$ and one of the following three cases holds:
\begin{enumerate}
\item[0)]  $u=x$, $v=y$ ($R\rightarrow S$ is unramified).
\item[1)]  $u=x$, $v=y^p\gamma+x\Sigma$ where $\gamma$ is a unit in $S$ and $\Sigma\in S$.
\item[2)] $u=\gamma x^p$, $v=y$ where $\gamma$ is a  unit in $S$.
\end{enumerate}

\end{Proposition}


Let $K\rightarrow L$ be an Artin-Schreier extension of two dimensional algebraic function fields over an algebraically closed field $k$ of characteristic $p>0$. Let $R\rightarrow S$ be an extension from a regular algebraic local ring of $K$ to a regular algebraic local ring of $L$ such that $S$ dominates $R$. 

Let $u, v$ be regular parameters in $R$ and $x,y$ be regular parameters in $S$. 
We will say that $R\rightarrow S$ is of type 0 with respect to these parameters if 
$$
\mbox{Type 0:}\,\,\,\,\, u=\gamma x, v=y\tau +x\Omega
$$
where $\gamma,\tau$ are units in $S$ and $\Omega\in S$, so that $R\rightarrow S$ is unramified.
We will say that $R\rightarrow S$ is of type 1 with respect to these parameters if 
$$
\mbox{Type 1:}\,\,\,\,\, u=\gamma x, v=y^p\tau +x\Omega
$$
where $\gamma,\tau$ are units in $S$ and $\Omega\in S$.
We will say that $R\rightarrow S$ is of type 2 with respect to these parameters if 
$$
\mbox{Type 2:}\,\,\,\,\, u=\gamma x^p, v=y\tau +x\Omega
$$
where $\gamma,\tau$ are units in $S$ and $\Omega\in S$.

These definitions are such that if one these types hold, and $\overline u,\overline v$ are regular parameters in $R$, $\overline x,\overline y$ are regular parameters in $S$  
such that $\overline u$ is a unit in $R$ times $u$ and $\overline x$ is a unit in $S$ times $x$ then $R\rightarrow S$ is of the same type for the new parameters $\overline u,\overline v$ and $\overline x,\overline y$.


In the construction of our example (Theorem \ref{Example3}), we will make use of some results from \cite{C3}.

\begin{Theorem}\label{TheoremA}(\cite[Theorem 4.1]{C3}) Suppose that $R\rightarrow S$  is of type 1 with respect to regular parameters $x, y$ in $S$ and $u,v$ in $R$ and that $J(S/R)=(x^{\overline c})$. Let $\overline x=u$, $\overline y=y-g(\overline x)$ 
where $g(\overline x)\in k[\overline x]$ is a polynomial with zero constant term, 
so that $\overline x,\overline y$ are regular parameters in $S$. Computing the Jacobian determinate $J(S/R)$, we see that 
$$
u=\overline x, v=\overline y^p\gamma+\overline x^{\overline c}\overline y\tau+f(\overline x)
$$
 where $\gamma,\tau$ are unit series in $\hat S$ and $f(\overline x)=\sum e_i\overline x^i\in k[[\overline x]]$.
Make the change of variables $\overline v=v-\sum e_iu^i$ where the sum is over $i$ such that $i\le \frac{pq}{m}$ so that $u,\overline v$ are regular parameters in $R$.

Suppose that $m,q$ are positive integers with $m>1$ and $\mbox{gcd}(m,q)=1$. Let $\alpha$ be a nonzero element of $k$.
Consider the sequence of quadratic transforms $S\rightarrow S_1$ so that $S_1$ has regular parameters $x_1,y_1$ defined by
$$
\overline x=x_1^{m}(y_1+\alpha)^{a'}, \overline y=x_1^{q}(y_1+\alpha)^{b'}
$$
where  $a',b'\in \NN$ are such that $mb'-q a'=1$. 

 We have that $R\rightarrow S$ is of type 1 with respect to the regular parameters $\overline x,\overline y$ and $u,v$.
Let $\sigma=\mbox{gcd}(m,p q)$ which is 1 or $p$. 

There exists a unique   sequence of quadratic transforms $R\rightarrow R_1$ such that $R_1$ has regular parameters $u_1,v_1$ defined by 
$$
u=u_1^{\overline m}(v_1+\beta)^{c'}, \overline v=u_1^{\overline q}(v_1+\beta)^{d'}
$$
with $0\ne\beta\in k$
giving
  a commutative diagram of homomorphisms
$$
\begin{array}{lll}
R_1&\rightarrow &S_1\\
\uparrow&&\uparrow\\
R&\rightarrow &S
\end{array}
$$
such that $R_1\rightarrow S_1$ is quasi finite. 
We  have that  $J(S_1/R_1)=(x_1^{c_1})$ for some positive integer $c_1$ and $R_1\rightarrow S_1$ is quasi finite. 
Further:

\begin{enumerate}
 \item[0)] If $\frac{q}{m}\ge  \frac{\overline c}{p-1}$  then $R_1\rightarrow S_1$ is of type 0.
\item[1)] If $\frac{q}{m}< \frac{\overline c}{p-1}$ and $\sigma=1$ then $R_1\rightarrow S_1$ is of type 1 and
$$
\left(\frac{c_{1}}{p-1}\right)=\left(\frac{\overline c}{p-1}\right)m-q.
$$  
\item[2)] If $\frac{q}{m}< \frac{\overline c}{p-1}$ and $\sigma=p$ then $R_1\rightarrow S_1$ is of type 2 and 
$$
\left(\frac{c_{1}}{p-1}\right)=\left(\frac{\overline c}{p-1}\right)m-q+1.
$$
\end{enumerate}
In cases 1) and 2), $m=\sigma \overline m$, $pq=\sigma\overline q$ and  $\overline mc'-\overline q d'=1$. 
\end{Theorem}

\begin{Theorem}\label{TheoremB}(\cite[Theorem 4.3]{C3}) Suppose that $R\rightarrow S$  is of type 2 with respect to regular parameters $x,y$ in $S$ and $u,v$ in $R$ and that $J(S/R)=(x^{\overline c})$.  Let $g(u)\in k[u]$ be a polynomial with no constant term. Make the change of variables, letting $\overline v=v-g(u)$ and $\overline y=\overline v$, so that $x,\overline y$ are regular parameters in $S$ and $u,\overline v$ are regular parameters in $R$.

Suppose that $m,q$ are positive integers with $\mbox{gcd}(m,q)=1$. Let $\alpha$ be a nonzero element of $k$.
Consider the sequence of quadratic transforms $S\rightarrow S_1$ so that $S_1$ has regular parameters $x_1,y_1$ defined by
$$
x=x_1^{m}(y_1+\alpha)^{a'}, \overline y=x_1^{q}(y_1+\alpha)^{b'}
$$
where  $a',b'\in \NN$ are such that $mb'-q a'=1$. 

Let $\sigma=\mbox{gcd}(pm,q)$ which is 1 or $p$. There exists a unique sequence of quadratic transforms $R\rightarrow R_1$ such that $R_1$ has regular parameters $u_1,v_1$ defined by 
$$
u=u_1^{\overline m}(v_1+\beta)^{c'}, \overline v=u_1^{\overline q}(v_1+\beta)^{d'}
$$
where $pm=\sigma\overline m$, $q=\sigma \overline q$, $\overline md'-c'\overline q=1$ and $0\ne\beta\in k$,
 giving
  a commutative diagram of homomorphisms
$$
\begin{array}{lll}
R_1&\rightarrow &S_1\\
\uparrow&&\uparrow\\
R&\rightarrow &S
\end{array}
$$
such that $R_1\rightarrow S_1$ is quasi finite. 
We have that  $J(S_1/R_1)=(x_1^{c_1})$ for some positive integer $c_1$.   Further: 
\begin{enumerate}
\item[1)] If  $\sigma=1$ then $R_1\rightarrow S_1$ is of type 1 and
$$
\left(\frac{c_{1}}{p-1}\right)=\left(\frac{\overline c}{p-1}\right)m-m.
$$  
\item[2)] If  $\sigma=p$ then $R_1\rightarrow S_1$ is of type 2 and 
$$
\left(\frac{c_{1}}{p-1}\right)=\left(\frac{\overline c}{p-1}\right)m-m+1.
$$
\end{enumerate}

\end{Theorem}

A proof of the following proposition is given in \cite[Proposition 7.9]{C3}. More general results are proven in \cite{KR}.
 
 \begin{Proposition}\label{Prop100} (Kuhlmann and Piltant, \cite{KP}) Suppose that $K$ and $L$ are two dimensional algebraic function fields over an algebraically closed field $k$ of characteristic $p>0$ and $K\rightarrow L$ is an Artin-Schreier extension. Let $\omega$ be a rational rank one nondiscrete valuation of $L$ and let $\nu$ be the restriction of $\omega$ to $K$. Suppose that $L$ is a defect extension of $K$.

Suppose that $R$ is a regular algebraic local ring of $K$ and $S$ is a regular algebraic local ring of $L$ such that $\omega$ dominates $S$, $S$ dominates $R$ and $R\rightarrow S$ is of type 1 or 2.
Inductively applying Theorems \ref{TheoremA} and  \ref{TheoremB}, we construct a diagram where the horizontal sequences are birational extensions of regular local rings
\begin{equation}\label{eq23*}
\begin{array}{ccccccc}
S=S_0&\rightarrow &S_1&\rightarrow & S_2&\rightarrow & \cdots\\
\uparrow&&\uparrow&&\uparrow&&\\
R=R_0&\rightarrow &R_1&\rightarrow & R_2&\rightarrow & \cdots
\end{array}
\end{equation}
with $\cup_{i=1}^{\infty}S_i=\mathcal O_{\omega}$. 
Further assume that for each map $R_i\rightarrow S_i$, there are regular parameters $u,v$ in $R_i$ and $x,y$ in $S_i$ such that one of the following forms hold:
\begin{equation}\label{eq20*}
u=x, v =f
\end{equation}
where $\dim_kS_i/(x,f)=p$, or
\begin{equation}\label{eq21*}
u=\delta x^p, v = y
\end{equation}
where $\delta$ is a unit in $S_i$ and  in both cases that  $x=0$ is a local equation of the critical locus of $\mbox{Spec}(S_i)\rightarrow \mbox{Spec}(R_i)$. 
Let
 $$
 J_i=J(S_i/R_i)=(\frac{\partial u}{\partial x}\frac{\partial v}{\partial y}-\frac{\partial u}{\partial y}\frac{\partial v}{\partial x})
 $$
 be the Jacobian ideal of the map $R_i\rightarrow S_i$.
 
 Then the distance $\mbox{dist}(\omega/\nu)$ is computed by the formula
 $$
 -{\rm dist}(\omega/\nu)=\frac{1}{p-1}\inf_i\{\omega(J(S_i/R_i))\}   
 $$
 where the infimum is over the $R_i\rightarrow S_i$ in the sequence  (\ref{eq23*}).
\end{Proposition}

\section{An example of a tower of independent defect extensions in which strong local monomialization doesn't hold}

\begin{Theorem}\label{Example3} There exists a tower $(K,\nu)\rightarrow (L,\omega)\rightarrow (M,\mu)$ of independent defect Artin-Schreier extensions of valued two dimensional algebraic function fields over an algebraically closed field $k$ of  characteristic $p>0$ such that there exist algebraic regular local rings $A$ of $K$ and $C$ of $M$ such that $\mu$ dominates $C$ and $C$ dominates $A$ but strong local monomialization along $\mu$ does not hold above $A\rightarrow C$. 
\end{Theorem} 

\begin{Remark}\label{RemarkA'} Let $\delta\in \RR_{\ge 0}$ be a fixed ratio. Suppose that $R\rightarrow S$ is of type 1. By taking $m$ and $q$ sufficiently large in Theorem \ref{TheoremA} such that   $R_1\rightarrow S_1$ is of type 2, we can achieve that $v_1=\lambda y_1+g(x_1)$ where $\lambda$ is a unit in $S_1$ and the order of $g(x_1)$ is arbitrarily large.  Suppose that $R\rightarrow S$ is of type 2. By taking $m$ and $q$ sufficiently large in Theorem \ref{TheoremB}  such that  $R_1\rightarrow S_1$ is of type 1 we can achieve that $v_1=y_1^p\gamma+x_1^{c_1}y_1\tau+f(x_1)$ where $\gamma$ and $\tau$ are unit series in $S_1$ and the order of $f(x_1)$ is arbitrarily large. In both cases, we can choose $m$ and $q$ so that $\frac{q}{m}$ is arbitrarily close to $\delta$.
\end{Remark}

\begin{Remark}\label{RemarkB'}
In Theorem \ref{TheoremB}, we have an expression $\overline v=\tau y+f(x)$ where $\tau$ is a unit in $S$. Suppose that $m$ and $q$ are positive integers with $\mbox{gcd}(m,q)=1$ and such that $\mbox{ord }f(x)>\frac{q}{m}$. Then the proof of Theorem \ref{TheoremB} extends to show that the conclusions of Theorem \ref{TheoremB} hold with $\overline y$ replaced with $y$.
\end{Remark}

We now give the proof of Theorem \ref{Example3}.
\begin{proof} 
Let $K$ be a two dimensional algebraic function field over an algebraically closed field, and let $R_{-2}$ be a two dimensional algebraic regular local ring of $K$. Let $u_{-2},v_{-2}$ be regular parameters in $R_{-2}$.


Let $e$ be a positive integer. Let $c_{-2}=(p-1)e$. Let $\Theta$ be a root of the Artin-Schreier polynomial $X^p-X-v_{-2}u_{-2}^{-pe}$. Let $L=K(\Theta)$. Set $x_{-2}=u_{-2}$, $y_{-2}=u_{-2}^e\Theta$. Let $S_{-2}=R_{-2}[y_{-2}]_{(x_{-2},y_{-2})}$, which is an algebraic regular local ring of $L$ which dominates $R_{-2}$. The regular parameters $x_{-2},y_{-2}$ in $S_{-2}$ satisfy $u_{-2}=x_{-2}, v_{-2}=y_{-2}^p-x_{-2}^{e(p-1)}y_{-2}$, so that the extension $R_{-2}\rightarrow S_{-2}$ is of type 1. We have that 
$J(S_{-2}/R_{-2})=(x_{-2}^{c_{-2}})$, with
$\frac{c_{-2}}{p-1}>0$.

We first construct a commutative diagram
$$
\begin{array}{lll}
S_{-2}&\rightarrow& S_{-1}\\
\uparrow&&\uparrow\\
R_{-2}&\rightarrow &R_{-1}
\end{array}
$$
using Theorem \ref{TheoremA} so that $R_{-1}\rightarrow S_{-1}$ is of type 2.
Let $\Sigma$ be a root of the Artin-Schreier polynomial $X^p-X-y_{-1}x_{-1}^{-pe}$. Let $M=L(\Sigma)$. Set $z_{-1}=x_{-1}$, $w_{-1}=x_{-1}^e\Sigma$. Let $T_{-1}=S_{-1}[w_{-1}]_{(z_{-1},w_{-1})}$, which is an algebraic regular local ring of $M$ which dominates $S_{-1}$. The regular parameters $z_{-1},w_{-1}$ in $T_{-1}$ satisfy $x_{-1}=z_{-1}, y_{-1}=w_{-1}^p-z_{-1}^{e(p-1)}w_{-1}$, so that the extension $S_{-1}\rightarrow T_{-1}$ is of type 1. We have that $J(T_{-1}/S_{-1})=(z_{-1}^{c_{-1}'})$, with
$\frac{c_{-1}'}{p-1}>0$.

From Theorems \ref{TheoremA} and \ref{TheoremB}, we construct
$$
\begin{array}{ccc}
T_{-1}&\rightarrow &T_0\\
\uparrow&&\uparrow\\
S_{-1}&\rightarrow &S_0\\
\uparrow&&\uparrow\\
R_{-1}&\rightarrow&R_0\\
\end{array}
$$
such that $R_0\rightarrow S_0$ is of type 1 and $S_0\rightarrow T_0$ is of type 2. Explicitely, $R_{-1},R_0, S_{-1}, S_0, T_{-1}, T_0$ have respective regular parameters $(u_{-1},v_{-1})$, $(u_0,v_0)$, $(x_{-1},y_{-1})$, $(x_0,y_0)$ and $(z_{-1}, w_{-1})$, $(z_0,w_0)$ which are related by equations
$$
\begin{array}{l}
u_{-1}=u_0^{pm_0}(v_0+\beta_0)^{d_0'}, v_{-1}=u_0^{q_0}(v_0+\beta_0)^{e_0'}\\
x_{-1}=x_0^{m_0}(y_0+\alpha_0)^{a_0'}, y_{-1}=x_0^{q_0}(w_0+\alpha_0)^{g_0'}\\
z_{-1}=z_0^{pm_0}(v_0+\gamma_0)^{f_0'}, w_{-1}=z_0^{q_0}(w_0+\gamma_0)^{g_0'}
\end{array}
$$
where $p\not| q_0$ and $\frac{q_0}{pm_0}<\frac{c_{-1}'}{p-1}$ where $J(T_{-1}/S_{-1})=(z_{-1}^{c_{-1}'})$.

By Remarks \ref{RemarkA'} and \ref{RemarkB'}, we can construct $R_0\rightarrow S_0\rightarrow T_0$ so that we have expressions $y_0=\lambda_0w_0+g_0(z_0)$ where $\lambda_0$ is a unit in $T_0$ and $\mbox{ord } g_0(z_0)$ is arbitrarily large and $v_0=\sigma_0y_0^p+\tau_0x_0^{c_0}y_0+f_0(x_0)$ where $\sigma_0,\tau_0$ are units in $S_0$ and $\mbox{ord } f_0(x_0)$ is arbitrarily large.

We will inductively construct a commutative diagram within $K\rightarrow L\rightarrow M$ of two dimensional regular algebraic local rings
\begin{equation}\label{eq5**}
\begin{array}{ccccccc}
T_0&\rightarrow &T_1&\rightarrow& T_2&\rightarrow &\cdots\\
\uparrow&&\uparrow&&\uparrow&&\\
S_0&\rightarrow&S_1&\rightarrow &S_2&\rightarrow &\cdots\\
\uparrow&&\uparrow&&\uparrow&&\\
R_0&\rightarrow&R_1&\rightarrow &R_2&\rightarrow&\cdots\\
\end{array}
\end{equation}
such that $R_i\rightarrow S_i$ is of type 1 if $i$ is even and is of type 2 if $i$ is odd, $S_i\rightarrow T_i$ is of type 2 if $i$ is even and is of type 1 if $i$ is odd. Further, valuations $\nu$, $\omega$ and $\mu$ of the respective function fields $K$, $L$ and $M$ determined by these sequences are such that $K\rightarrow L$ and $L\rightarrow M$ are independent defect extensions. 
We will have that $R_i$ has regular parameters $(u_i,v_i)$, $S_i$ has regular parameters $(x_i,y_i)$ and $T_i$ has regular parameters $(z_i,w_i)$ such that
$$
u_i=u_{i+1}^{\overline m_{i+1}}(v_{i+1}+\beta_{i+1})^{d_{i+1}'}, v_i=u_{i+1}^{\overline q_{i+1}}(v_{i+1}+\beta_{i+1})^{e_{i+1}'},
$$
$$
x_i=x_{i+1}^{m_{i+1}}(y_{i+1}+\alpha_{i+1})^{a_{i+1}'}, y_i=x_{i+1}^{q_{i+1}}(y_{i+1}+\alpha_{i+1})^{b_{i+1}'},
$$
$$
z_i=z_{i+1}^{m_{i+1}'}(w_{i+1}+\gamma_{i+1})^{f_{i+1}'}, w_i=z_{i+1}^{q_{i+1}'}(w_i+\gamma_{i+1})^{g_{i+1}'}
$$
with $\overline m_i, m_i$ and $m_i'$ larger than 1 for all $i$.

Let $J(S_i/R_i)=(x_i^{c_i})$ and $J(T_i/S_i)=(z_i^{c_i'})$.

If $i$ is even, then 
$m_{i+1}=p\overline m_{i+1}, m_{i+1}'=\overline m_{i+1},q_{i+1}=\overline q_{i+1}, q_{i+1}'=q_{i+1}$ and 
$$
\frac{q_{i+1}}{m_{i+1}}<\frac{c_i}{p-1}.
$$

If $i$ is odd, then $\overline m_{i+1}=pm_{i+1}, m_{i+1}'=\overline m_{i+1}, q_{i+1}=\overline q_{i+1}, q_{i+1}'=q_{i+1}$ and 
$$
\frac{q_{i+1}'}{m_{i+1}'}<\frac{c_i'}{p-1}.
$$

In our construction, if $r$ is even, we will have that
\begin{equation}\label{eq100*}
y_r=\lambda_rw_r+g_r(z_r)
\end{equation}
where $\lambda_r$ is a unit in $T_r$ and $\mbox{ord } g_r(z_r)$ is arbitrarily large and
\begin{equation}\label{eq101}
v_r=\sigma_ry_r^p+\tau_rx_r^{c_r}y_r+f_r(x_r)
\end{equation}
where $\sigma_r,\tau_r$ are units in $S_r$ and $\mbox{ord }f_r(x_r)$ is arbitrarily large. 
If $r$ is even, we will have
\begin{equation}\label{eq106}
y_r=\sigma_rw_r^p+\tau_rz_r^{c_r'}w_r+f(z_r)
\end{equation}
where $\sigma_r,\tau_r$ are units in $T_r$ and $\mbox{ord }f(z_r)$ is arbitrarily large and
\begin{equation}\label{eq107}
v_r=\lambda_ry_r+g_r(x_r)
\end{equation}
where $\lambda_r$ is a unit in $S_r$ and $\mbox{ord }g_r(x_r)$ is arbitrarily large.

Suppose that $r$ is even, and we have constructed $R_r\rightarrow S_r\rightarrow T_r$. We will construct
$$
\begin{array}{ccccc}
T_r&\rightarrow&T_{r+1}&\rightarrow&T_{r+2}\\
\uparrow&&\uparrow&&\uparrow\\
S_r&\rightarrow&S_{r+1}&\rightarrow&S_{r+2}\\
\uparrow&&\uparrow&&\uparrow\\
R_r&\rightarrow&R_{r+1}&\rightarrow &R_{r+2}.
\end{array}
$$
There exists an integer $\lambda(r+1)>1$ and $q_{r+1}\in \ZZ_+$ such that $\mbox{gcd}(q_{r+1},p)=1$ and
\begin{equation}\label{eq64*}
\frac{c_r}{p-1}
>\frac{q_{r+1}}{p^{\lambda(r+1)}}>\frac{c_r}{p-1}-\frac{1}{2^{r+1}}m_1\cdots m_{r}.
\end{equation}
In fact, we can find $\lambda(r+1)$ arbitrarily large satisfying the inequality.
Set $m_{r+1}=p^{\lambda(r+1)}$. 
We have that $\frac{q_{r+1}}{m_{r+1}}<\frac{c_r}{p-1}$ with $\mbox{gcd}(m_{r+1},pq_{r+1})=p$. 
This choice of $m_{r+1}$ and $q_{r+1}$ (along with a choice of $0\ne\alpha_{r+1}\in k$) determines $S_r\rightarrow S_{r+1}$. We have an expression
$v_r=\sigma_r y_r^p+\tau_r x_r^{c_r}y_r+f_r(x_r)$ where $\mbox{ord }f_r(x_r)$ is arbitrarily large. In particular, we can assume that
$\mbox{ord }f_r(x_r)>\frac{pq_{r+1}}{m_{r+1}}$. Then $R_r\rightarrow R_{r+1}$ is defined as desired by Theorem \ref{TheoremA}.
By Remark \ref{RemarkA'}, since we can take $\lambda(r+1)$ to be arbitrarily large, we can assume that $v_{r+1}=\lambda_{r+1}y_{r+1}+g_{r+1}(x_{r+1})$ where $\mbox{ord }g_{r+1}(x_{r+1})$ is arbitrarily large. 

By Remark \ref{RemarkB'} and Theorem \ref{TheoremB}, $T_r\rightarrow T_{r+1}$ is defined as desired, with $m_{r+1}'=\frac{m_{r+1}}{p}$, $q_{r+1}'=q_{r+1}$. Since we can take $\lambda(r+1)$ to be arbitrarily large, we can assume that
$y_{r+1}=\sigma_{r+1}w_{r+1}^p+\tau_{r+1}z_{r+1}^{c_{r+1}'}w_r+f_{r+1}(z_{r+1})$ where $\mbox{ord }f_{r+1}(z_{r+1})$ is arbitrarily large. 

We have defined a commutative diagram
\begin{equation}\label{eq103}
\begin{array}{ccc}
T_r&\rightarrow &T_{r+1}\\
\uparrow&&\uparrow\\
S_r&\rightarrow &S_{r+1}\\
\uparrow&&\uparrow\\
R_r&\rightarrow & R_{r+1}
\end{array}
\end{equation}
with the desired properties; in particular, $R_{r+1}\rightarrow S_{r+1}$ is of type 2 with
$$
\frac{c_{r+1}}{p-1}=\left(\frac{c_r}{p-1}\right)m_{r+1}-q_{r+1}+1
$$
and $S_{r+1}\rightarrow T_{r+1}$ is of type 1, with 
$$
\frac{c_{r+1}'}{p-1}=\frac{c_r'}{p-1}m_{r+1}'-m_{r+1}'.
$$


Now choose $q_{r+2}'$, $m_{r+2}'=p^{\lambda(r+2)}$ such that $p\not|q_{r+2}'$ and
\begin{equation}\label{eq1**}
\frac{c_{r+1}'}{p-1}>\frac{q_{r+2}'}{m_{r+2}'}>\frac{c_{r+1}'}{p-1}-\frac{1}{2^{r+2}}m_1'\cdots m_{r+1}'.
\end{equation} 
We can take $\lambda(r+2)$ arbitrarily large.
Set $m_{r+2}=\frac{m_{r+2}'}{p}=p^{\lambda(r+2)-1}$, $q_{r+2}=q_{r+2}'$. By (\ref{eq1**}),
$\frac{q_{r+2}'}{m_{r+2}'}<\frac{c_{r+1}'}{p-1}$.

Now construct, as in the construction of (\ref{eq103}), using Theorems \ref{TheoremA} and \ref{TheoremB} and Remark \ref{RemarkB'} and these values of $m_{r+2}$ and $q_{r+2}$,
$$
\begin{array}{ccc}
T_{r+1}&\rightarrow &T_{r+2}\\
\uparrow&&\uparrow\\
S_{r+1}&\rightarrow &S_{r+2}\\
\uparrow&&\uparrow\\
R_{r+1}&\rightarrow & R_{r+2},
\end{array}
$$
so that $R_{r+2}\rightarrow S_{r+2}$ is of type 1 and $S_{r+2}\rightarrow T_{r+2}$ is of type 2. By Remark \ref{RemarkA'}, we obtain expressions (\ref{eq100*}) and (\ref{eq101}) for $r+2$.

 By induction, we construct the diagram (\ref{eq5**}).

  Let $A=R_0$ and $C=T_0$. We will show that strong local monomialization doesn't hold above $A\rightarrow C$ along $\mu$. Suppose that $R'\rightarrow T'$ has a strongly monomial form above $A\rightarrow C$. Then $R'$ has regular parameters $u',v'$ and $T'$ has regular parameters $z',w'$ such that $u'=\lambda (z')^m$ and $v'=w'$ where $m\in \ZZ_{>0}$ and $\lambda$ is a unit in $T'$. We will show that this cannot occur. There exists a commutative diagram
  $$
  \begin{array}{ccccc}
  T_s&\rightarrow&T'&\rightarrow&T_{s+1}\\
  \uparrow&&\uparrow&&\uparrow\\
  R_s&\rightarrow&R'&\rightarrow&R_{s+1}
  \end{array}
  $$
  for some $s$. The ring $T'$ has regular parameters $\overline z,\overline w$ such that 
  \begin{equation}\label{eq108}
  z_s=\overline z^a\overline w^b, w_s=\overline z^c\overline w^d
  \end{equation}
   for some $a,b,c,d\in \NN$ with $ad-bc=\pm1$, and $R'$ has regular parameters $\overline u,\overline v$ such that $u_s=\overline u^{\overline a}\overline v^{\overline b}$, $v_s=\overline u^{\overline c}\overline v^{\overline d}$, where $\overline a \overline d-\overline b \overline c=\pm 1$. We have an expression 
  \begin{equation}\label{eq102}
  u_s=\alpha z_s^p, v_s=\beta w_s^p+\Omega
  \end{equation}
  where $\alpha,\beta$ are units in $T_s$ and where 
  \begin{equation}\label{eq109}
  \Omega=\epsilon z_s^{pc_s}w_s+M
  \end{equation}
  or
  \begin{equation}\label{eq110}
  \Omega=\epsilon z_s^{c_s'}w_s+M
  \end{equation}
  where $\epsilon\in T_s$ is a unit and $M$ is a sum of monomials in $z_s,w_s$ of high order in $z_s$. Further,
  $\mu(w_s^p)<\mu(z_s^{pc_s}w_s)$ in (\ref{eq109}) and
  $\mu(w_s^p)<\mu(z_s^{c_s'}w_s)$ in (\ref{eq110}).
  
  In particular, $R_s\rightarrow T_s$ is not a strongly monomial form. 
  
   Substituting (\ref{eq108}) into $u_s$ and $v_s$ in (\ref{eq102}), we have 
  \begin{equation}\label{eq113}
  u_s=\alpha \overline z^{ap}\overline w^{bp}, v_s=\beta \overline z^{cp}\overline w^{dp}+\Omega.
  \end{equation}
  We necessarily have that $u_s|v_s$ or $v_s|u_s$ in $T'$. 
  
  First suppose that $c\ge a$ and $d\ge b$. Then we have that
  $$
  u_s=\alpha \overline z^{ap}\overline w^{bp},
  \frac{v_s}{u_s}=\beta\overline z^{(c-a)p}\overline w^{(d-b)p}+\frac{\Omega}{\alpha \overline z^{ap}\overline w^{bp}}
  $$
  giving an expression of the form (\ref{eq113}).
  We will show that this is not a strongly monomial form. If it is, then we must have that $a=0$ or $b=0$ so that either 
  \begin{equation}\label{eq111}
  z_s=\overline w, w_s=\overline z \overline w^d
  \end{equation}
  or
  \begin{equation}\label{eq112}
  z_s=\overline z, w_s=\overline z^c \overline w
  \end{equation}
  and we must have that $\frac{\Omega}{u_s}$ is part of a regular system of parameters in $T'$. Substituting into (\ref{eq109}) or (\ref{eq110}), we see that this  cannot occur except possibly in the case that (\ref{eq110}) holds and $\frac{z_s^{c_s'}w_s}{u_s}$ is part of a regular system of parameters in $T'$.

 Suppose that  (\ref{eq110}) and (\ref{eq111}) hold with
  $$
  \frac{z_s^{c_s'}w_s}{u_s}=\frac{\overline w^{c_s'+d}\overline z}{\alpha \overline w^p}
  $$
  being part of a regular system of parameters in $T'$. Now in this case, $\mu(w_s)>\mu(z_s)$ and $\mu(w_s^p)<\mu(z_s^{c_s'}w_s)$ so $p\le c_s'$. Thus $\frac{\overline w^{c_s'+d}\overline z}{\alpha \overline w^p}$ cannot be part of a regular system of parameters in $T'$.   A similar argument shows that we do not obtain a strongly monomial form when (\ref{eq110}) and (\ref{eq112}) hold. 
  
  Suppose that $c<a$ and $d<b$. Then we have expressions
  $$
  v_s=\gamma\overline z^{cp}\overline w^{dp},
  \frac{u_s}{v_s}=\alpha \gamma^{-1}\overline z^{(a-c)p}\overline w^{(b-d)p}
  $$
  where $\gamma\in T'$ is a unit, giving an expression of the form of (\ref{eq113}), which is not strongly monomial. Thus we reduce to the case where $(c-a)(d-b)<0$. We then have that $u_s\not | v_s$ since $u_s\not |\overline z^{cp}\overline w^{dp}$. Suppose that $v_s|u_s$. Then $v_s=\lambda \overline z^{cp}\overline w^{dp}$ where $\lambda$ is a unit in $T'$. But this is impossible since $(c-a)(d-b)<0$.
 Thus $R'\rightarrow T'$ has a form (\ref{eq113}) with $a,b,c,d>0$ and so 
 cannot be a strongly monomial form. We have established that strong local monomialization along $\mu$ does not hold above $A\rightarrow C$.

From Theorem \ref{TheoremA}, we have that 
\begin{equation}\label{eq60*}
\left(\frac{c_{r+1}}{p-1}\right)\frac{1}{m_1\cdots m_{r+1}}=\left(\frac{c_r}{p-1}\right)\frac{1}{m_1\cdots m_{r}}
-\frac{q_{r+1}}{m_{r+1}}\left(\frac{1}{m_1\cdots m_{r}}\right)+\frac{1}{m_1\cdots m_{r+1}}.
\end{equation}
Then from Theorem \ref{TheoremB}, we have that
$$
\frac{c_{r+2}}{p-1}=\left(\frac{c_{r+1}}{p-1}\right)m_{r+2}-m_{r+2},
$$
and so
\begin{equation}\label{eq61*}
\left(\frac{c_{r+2}}{p-1}\right)\frac{1}{m_1\cdots m_{r+2}}
=\left(\frac{c_r}{p-1}\right)\frac{1}{m_1\cdots m_{r}}-\frac{q_{r+1}}{m_1\cdots m_{r+1}}.
\end{equation}
By equation (\ref{eq64*})  we have
\begin{equation}\label{eq65*}
\frac{1}{2^{r+1}}>\left(\frac{c_r}{p-1}\right)\frac{1}{m_1\cdots m_{r}} -\left(\frac{q_{r+1}}{m_{r+1}}\right)\frac{1}{m_1\cdots m_{r}}>0.
\end{equation}


By Theorem \ref{TheoremA},
$$
\left(\frac{c_{r+2}'}{p-1}\right)\frac{1}{m_1'\cdots m_{r+2}'}=\left(\frac{c_{r+1}'}{p-1}\right)\frac{1}{m_1'\cdots m_{r+1}'}
-\frac{q_{r+2}'}{m_1'\cdots m_{r+2}'}+\frac{1}{m_1'\cdots m_{r+2}'}
$$
and by Theorem \ref{TheoremB},
$$
\frac{c_{r+3}'}{p-1}=\left(\frac{c_{r+2}'}{p-1}\right)m_{r+3}'-m_{r+3}'.
$$
We thus have that
\begin{equation}\label{eq2**}
\left(\frac{c_{r+3}'}{p-1}\right)\frac{1}{m_1'\cdots m_{r+3}'}=\left(\frac{c_{r+1}'}{p-1}\right)\frac{1}{m_1'\cdots m_{r+1}'}-
\frac{q_{r+2}'}{m_1'\cdots m_{r+2}'}.
\end{equation}
Equation (\ref{eq1**}) implies
\begin{equation}\label{eq3**}
\frac{1}{2^{r+2}}>\left(\frac{c_{r+1}'}{p-1}\right)\frac{1}{m_1'\cdots m_{r+1}'}-\frac{q_{r+2}'}{m_1'\cdots m_{r+2}'}>0.
\end{equation}

Now $J(S_i/R_i)=(x_i^{c_i})$ and $x_0=x_i^{m_1\cdots m_i}$ so
$\omega(J(S_i/R_i))=\frac{c_i}{m_1\cdots m_i}\omega(x_0)$ and thus
by Proposition \ref{Prop100}, (\ref{eq61*}) and (\ref{eq65*}), we have that
$$
-{\rm dist}(\omega/\nu)=\frac{1}{p-1}\inf_i\{\omega(J(S_i/R_i))\}=0.
$$

We  have that 
$J(T_i/S_i)=(z_i^{c_i'})$ and $z_0=z_i^{m_1'\cdots m_i'}$ so
$\omega(J(T_i/S_i))=\frac{c_i'}{m_1'\cdots m_i'}\omega(z_0)$ and thus
by Proposition \ref{Prop100}, (\ref{eq2**}) and (\ref{eq3**}), we have that
$$
-{\rm dist}(\mu/\omega)=\frac{1}{p-1}\inf_i\{\omega(J(T_i/S_i))\}=0.
$$

\end{proof}

\end{document}